\documentclass[notimes]{georeport}

\usepackage{natbib}
\usepackage{graphicx}
\usepackage{color}
\usepackage{listings}
\usepackage{amsmath}
\usepackage{amsthm}
\usepackage{amssymb}
\usepackage{amsbsy}
\usepackage{float}
\usepackage{wasysym}
\usepackage{setspace}
\usepackage{bm}
\usepackage{algorithm}
\usepackage{algpseudocode}

      \newcommand{\lW}{W_{\lambda}}

      \newcommand{\lF}{F_{\lambda}}
      \newcommand{\lerr}{e_{\lambda}}
      \newcommand{\Ja}{J_{\alpha}}
      \newcommand{\tJa}{\tilde{J}_{\alpha}}	
      \newcommand{\bx}{{\bf x}}
      \newcommand{\bR}{{\bf R}}
      \newcommand{\aw}{w_{\alpha}}

      \newtheorem{theorem}{Theorem}
      \newtheorem{corollary}{Corollary}
      
      \newtheorem{proposition}{Proposition}

\setfigdir{Fig}

\begin{document}
\title{Error Bounds for Extended Source Inversion applied to an Acoustic Transmission Inverse Problem}
\author{William W. Symes\\
PO Box 43, Orcas, WA 98280 USA, {\tt symes@rice.edu}}

\lefthead{Symes}

\righthead{1D Acoustic Inversion: Theory}

\maketitle
\begin{abstract}
  A simple inverse problem for the wave equation requires
  determination of both the wave velocity in a homogenous acoustic
  material and the transient waveform of an isotropic point radiator,
  given the time history of the wavefield at a remote point in
  space. The duration (support) of the source waveform and the
  source-to-receiver distance are assumed known. A least squares
  formulation of this problem exhibits the ``cycle-skipping''
  behaviour observed in field scale problems of this type, with many
  local minima differing greatly from the global minimizer. An
  extended formulation, dropping the support constraint on the source
  waveform in favor of a weighted quadratic penalty, eliminates this
  misbehaviour. With proper choice of the weight operator, the
  velocity component at {\em any} local minimizer of this extended
  objective function differs from the global
  minimizer of the least-squares formulation by less than a linear
  combination of the source waveform support radius and data
  noise-to-signal ratio.

\end{abstract}

\section{Introduction}
Inverse problems based on wave equations (acoustic, elastic,
Maxwell's...) are commonplace in geophysics, nondestructive materials
testing, medical imaging, and other areas of science and engineering
in which wave motion plays an important role. These problems may be
formulated as nonlinear least squares problems, requiring optimization
of mean square misfit between observed and predicted data. Iterative
local numerical optimization methods (relatives of Newton's method)
are feasible solution approaches for the computationally large
variants of these nonlinear least squares problems that occur in
seismology and medical imaging. However, local optimization produces
only approximate stationary points, and these exhibit a tendency to
stagnate at physically irrelevant solutions, under typical conditions
of data acquisition. This ``cycle-skipping''
pathology can sometimes be overcome through use of information other than
measured wave data, however it remains a serious obstacle to effective
solution of field-scale inverse problems in wave propagation
\cite[]{GauTarVir:86,VirieuxOperto:09,Fichtner:10,Plessix:10,Schuster:17}.

Since numerical feasibility limits practical optimization approaches
to search for stationary points, it is natural to explore alternatives
to nonlinear least squares formulations for which approximate
stationary points might perforce be acceptable solutions. Many such
alternatives have been suggested and applied in both synthetic and
field data tests, in some cases with excellent results (see
\cite{HuangNammourSymesDollizal:SEG19},
\cite{PladysBrossierLiMetivier:GEO21} for recent reviews). However
numerical tests are inevitably ``a look in the rear view mirror'',
with no guarantee that the next example will not behave quite
differently. See \cite{Symes:2020} for an example of such an unpleasant
surprise. Theoretical guarantees that stationary points are
satisfactory solutions of specific inverse problems, in some specific
sense and under specified conditions, are almost entirely absent from
the literature.

This paper provides such a theoretical guarantee for a modification of
least-squares data fitting via {\em modeling operator extension},
applied to a very simple example: recovery of the (spatially
homogeneous) wave speed of an acoustic material model, together with
the waveform of an isotropic point radiator (``wavelet'') with
specified support, from the time history (``trace'') of the resulting
pressure field at a single remote location. This is perhaps the
simplest inverse wave propagation problem that exhibits
cycle-skipping. While far too simple to have immediate
practical application, it is a special case, or subproblem, of many
field-scale inverse problems in industrial and academic seismology. It
is also routinely used to illustrate cycle-skipping (see for example
\cite{VirieuxOperto:09}, Figure 7). Modeling operator extension ideas underlie
widely used techniques in seismic data processing, and have been
applied to inversion for several decades
\cite[]{geoprosp:2008}. The particular approach used in this work
is an example of {\em source extension}.
\cite{HuangNammourSymesDollizal:SEG19} review recent source extension
contributions, some of which are closely related to the approach to
single-trace transmission inversion studied here.

Model extension is only one of several ideas proposed to remedy
cycle-skipping. For instance, measuring data misfit with versions of the Wasserstein
metric from optimal transport theory, rather than the $L^2$ norm, also
shows promise in mitigating cycle-skipping
\cite[]{Metivier:GEO18,EngquistYang:GEO18,Ramos-Martinez:SEG18,Wang:SEG19,EngquistYang:CPAM21}. \cite{MahankaliYang:21}
have established conditions under which the Wasserstein misfit is
convex in a few parameters of simple models. Aside from the present
paper, their work is one of the very few to give a complete mathematical
justification (albeit also in very
restricted circumstances) for application of local optimization to a nonlinear
least squares alternative.

The next section precisely defines the single-trace transmission
inverse problem, and the following section confirms the
cycle-skipping behaviour of the obvious nonlinear least squares
formulation, that is, that stationary points exist at arbitrarily
large distance from any global minimizer (Theorem \ref{thm:fwi}). The
modeling operator extension used here simply consists in dropping the
support constraint, thus allowing fit to any data with any wavespeed,
by suitable choice of wavelet. The link between wave speed and data is
re-established by supplementing the mean-square data misfit function
with a suitable quadratic penalty on nonzero values of the wavelet far
from the specified support. I show that any stationary point of the
resulting penalty function deviates from the global solution of the
original, constrained nonlinear least squares problem by an amount
proportional to a combination of the support diameter and the relative
data error, and this bound is sharp
(Theorems \ref{thm:rampreallygood} and \ref{thm:mnoiseres}). A
solution of the original inverse problem (with wavelet support
constraint) is recovered by truncating the wavelet obtained in the
course of minimizing the penalty function (Theorem
\ref{thm:ipnoisesuf}).

The elementary
analysis developed here leans on the special features of the problem
under study to achieve sharper results than will likely be possible in
more general contexts. The Appendix recasts some of these
constructions in a form that applies, at least formally, to more
prototypical, larger-scale problems
\cite[]{Symes:IPTA14,tenKroode:IPTA14}, and suggests a route to
achieving a similar, mathematically complete understanding of
extension-based approaches to applied wave propagation inverse
problems.

\section{An Acoustic transmission Inverse Problem}
According to linear acoustics, the causal pressure field $p(\bx,t)$ due to an
isotropic point radiator at $\bx_s \in \bR^3$ with time-varying
intensity (``wavelet'') $w(t)$, propagating in an elastic fluid with
slowness (reciprocal velocity) $m$ and density $1$ (in appropriate
units), is the solution of the the initial value problem for the wave equation \cite[]{Frie:58}:
\begin{eqnarray}
  \label{eqn:awe}
  \left(m^2\frac{\partial^2 p}{\partial t^2} - \nabla^2\right) p(\bx,t) &=&
                                                                         w(t)\delta(\bx-\bx_s) \nonumber\\
  p(\bx,t)&=&0, t\ll 0.
\end{eqnarray}
The solution is well-known, see for instance
\cite{CourHil:62}, Chapter VI, section 12, equation 47:
\begin{equation}
  \label{eqn:homsol}
  p(\bx,t) = \frac{1}{4\pi |\bx-\bx_s|}w\left(t-m|\bx-\bx_s|\right),
\end{equation}
This trace of $p$ at $\bx_r \in \bR^3$ can be viewed as the result of
applying an m-dependent linear operator $F[m]$ to the wavelet $w$:
\begin{equation}
\label{eqn:mod}
F[m]w(t)  = p(\bx_r,t) = \frac{1}{4\pi r}w\left(t-mr\right) 
\end{equation}
The slowness $m$ must be positive, as follows from basic acoustics,
and in fact reside in a range characteristic of the
material model: for crustal rock, a reasonable choice would be
$m_{\rm min}=0.125, m_{\rm max}=0.6$ s/km. The pressure trace
$p(\bx_r,\cdot)$ is square-integrable if and only if the same is true
of the wavelet $w$. Since the square-integral of the pressure trace is
proportional to the accumulated energy transferred from the fluid to
the sensor \cite[]{SantosaSymes:00}, assume that $w \in L^2(\bR)$.

Natural choices for domain and
range of $F$ are thus
\begin{itemize}
\item $M=(m_{\rm min}, m_{\rm max}),\,0 < m_{\rm min} \le m_{\rm
    max}$;
\item $W = L^2(\bR)$;
\item $D=L^2([t_{\rm min},t_{\rm max}]),\, t_{\rm min}<t_{\rm max}$;
\item $F: M \times W \rightarrow D$ as specified in \ref{eqn:mod}.
\end{itemize}
It is immediately evident from these choices and from the definition
\ref{eqn:mod} that
\begin{equation}
  \label{eqn:mapprop}
  \mbox{for }m \in M, F[m] \mbox{ is bounded, and }\|F[m]\| =
  \frac{1}{4\pi r}.
\end{equation}

Since
all possible data lie in the range of $F[m]$ for any $m \in M$, some
restriction of the domain of $F$ is necessary in order that fitting
the data constrain $m$. The constraint employed in this work is the specification
of a maxium support radius $\lambda_{\rm max} >0$ (see the companion
paper \cite{SymesChenMinkoff:21} for a justification of this choice).

For
$\lambda \in (0,\lambda_{\rm max}]$, define
\begin{itemize}
\item $\lW = \{w \in W:
  \mbox{supp }w \subset [-\lambda,\lambda]\}$;
\item $\lF = F|_{M \times \lW}$.
\end{itemize}

In terms of this infrastructure, the inverse problem studied in 
this paper may be stated as

\begin{quote}
\noindent {\bf Inverse Problem:}
  given data $d \in D$, relative error level $\epsilon \in
  [0,1)$, and support radius $\lambda \in (0, \lambda_{\rm
    max}]$, find $(m,w) \in M \times \lW$ for which 
\begin{equation}
  \label{eqn:probstat0}  \|\lF[m]w-d\| \le \epsilon\|d\|,
\end{equation}
\end{quote}

Define the relative mean-square error $\lerr: M \times \lW \times D
\rightarrow \bR^+$ by
\begin{equation}
  \label{eqn:redms}
  \lerr[m,w;d]=\frac{1}{2}\|\lF[m]w-d\|^2/\|d\|^2,
\end{equation}
so that inequality \ref{eqn:probstat0} is equivalent to
\begin{equation}
  \label{eqn:probstat1}
  \lerr[m,w;d] \le \frac{1}{2}\epsilon^2,
\end{equation}
Minimization of $\lerr$ is the standard least-squares formulation of
the Inverse Problem defined above: if the global minimum value of
$\lerr$ is less than $\frac{1}{2}\epsilon^2$, then the global
minimizer is a solution of the Inverse Problem.

\noindent{\bf Remark:} The constraint $\epsilon < 1$ imposed on the
target noise level eliminates the obvious choice $(m,0)$, which
satisfies the data misfit constraint for any $m \in M$ if $\epsilon
\ge 1$. 

\noindent{\bf Remark:} I shall refer to the minimization of $\lerr$ as
``Full Waveform Inversion'' or ``FWI'', as this is the
terminology used in the seismology literature to identify this and
similar optimization problems.

The best case for data fitting
is clearly the one in which the data can be fit precisely: that is,
there exists $(m_*,w_*) \in M\times \lW$ so that
\begin{equation}
  \label{eqn:defdatanonoise}
  d=\lF[m_*]w_*.
\end{equation}
Such data $d$ is {\em noise-free}, in the range of the map $\lF$. For
such data a solution of the Problem Statement \ref{eqn:probstat0}
exists with arbitrarily small $\epsilon>0$.

\section{Full Waveform Inversion}
While $F$ is surjective, as noted above, it is
very far from injective. On the other hand, under a constraint that
will be assumed throughout, $\lF[m]$ is injective for each $m \in M$ (in fact, $4 \pi r
\lF[m]$ is an isometry):
\begin{proposition}
  \label{thm:fullrec}
  Suppose that 
  \begin{equation}
    \label{eqn:fullrec}
    [ m_{\rm min}r-\lambda_{\rm max}, m_{\rm max}r+\lambda_{\rm max}]
    \subset [t_{\rm min},t_{\rm max}].
  \end{equation}
  Then $\lF[m]$ is coercive for every $m \in M, \lambda \in
  (0,\lambda_{\rm max}]$.
\end{proposition}

\noindent{\bf Remark:} A useful consequence of the condition 
\ref{eqn:fullrec}: for every $m \in M$, 
\begin{equation}
  [-\lambda_{\rm max}, -\lambda_{\rm max}] \subset[ t_{\rm min}-mr , 
  t_{\rm max}-mr]. 
  \label{eqn:zeroinc}
\end{equation}

The first main result establishes the existence of large (100 \%)
residual local minimizers for the basic FWI objective $\lerr$, even
for noise-free data.
\begin{theorem}
  \label{thm:fwi}
  Suppose that $0 <\lambda\le \lambda_{\rm max}$,  $m_* \in M, w_*
  \in \lW, d=\lF[m_*]w_*$ is noise-free data per definition \ref{eqn:defdatanonoise},
  Under assumption \ref{eqn:fullrec}, for any $m \in M$ with $|m-m_*|r>2\lambda$,
\begin{equation}
  \label{eqn:isovpm}
 \min_w \lerr[m,w;d]=\lerr[m,0;d] = \frac{1}{2},
\end{equation}
and any such $(m,0)$ is a local minimizer of $\lerr$ with relative RMS
error = 1.0.
\end{theorem}

\begin{proof} From the definition \ref{eqn:mod},
\[
 \lerr[m,w;d] =  \frac{1}{32\pi^2
    r^2\|d\|^2}\int\,dt\,\left|w\left(t-mr\right)-w_*\left(t-m_*r\right)\right|^2
\]
Since $w_*, w$ vanish for $|t|>\lambda$,
$\lF[m_*]w_*(t)$ vanishes if $|t-m_*r|>\lambda$ and $\lF[m]w$ vanishes if $|t-mr|>\lambda$. So if $|mr-m_*r|
= |m-m_*|r > 2\lambda$, then $|t-mr|+|t-m_*r| \ge |mr-m_*r| >
2\lambda$ so either $|t-mr|>\lambda$ or $|t-m_*r|>\lambda$, that is,
either $\lF[m]w(t)=0$ or $\lF[m_*]w_*=0$. Therefore $\lF[m]w$ and
$\lF[m_*]w_*$ are orthogonal in the sense of the $L^2$ inner product
$\langle \cdot,\cdot \rangle_D$ on $D$:
\begin{equation}
  \label{eqn:ortho}
  |m- m_*|r > 2\lambda \,\, \Rightarrow \,\, \langle F[m]w,
  F[m_*]w_*\rangle_D = 0
\end{equation}
But $d = \lF[m_*]w_*$, so this is the same as saying that $d$ is
orthogonal to $F[m]w$. So conclude after a minor manipulation that
\[
  |m- m_*|r > 2\lambda \,\, \Rightarrow \,\, \lerr[m,w;d]=\frac{1}{32\pi^2 
    r^2\|d\|^2}(\|w\|^2 + \|w_*\|^2).
\]
\begin{equation}
  \label{eqn:iso}
  = \frac{1}{2}\left(\frac{\|w\|^2}{\|w_*\|^2} + 1 \right)
\end{equation}
That is, for slowness $m$ in error by more than $2\lambda/r$ from the 
target slowness $m_*$, the means square error (FWI objective) $\lerr$ is independent of
$m$, and its minimum over $w$ is attained for $w=0$
\end{proof}

Therefore local minimizers of $\lerr$ abound, as far as you like from the
global minimizer $(m_*,w_*)$. Local exploration of the FWI objective
$e$ gives no useful information whatever about constructive search
directions, and descent-based optimization tends to fail if the
initial estimate $m_0$ is in error by more than $2\lambda/r$. In fact the actual behaviour of FWI itererations is worse
(failure if $m_0$ is in error by ``half a wavelength''), as follows
from a more refined analysis of the cycle-skipping local behaviour of $\lerr$ near its
global minimizer.

\section{Extended Source Inversion}
The phenomenon explained in the last section can be avoided by
reformulating the inverse problem via an extended modeling operator
and a soft (penalty) constraint to replace the support
requirement. The extension simply amounts to dropping the support
constraint, and replacing $\lW$ and $\lF$ by $W$ and $F$.

The hard
support constraint implicit in the choice of $\lW$ as domain for the
modeling operator is replaced by a soft constraint in the form of a
quadratic penalty, with weight operator $A:W \rightarrow W$.
The choices for the penalty operator $A$ considered here are scalar 
multiplication operators on $W$ defined by a choice of multiplier $a \in L^{\infty}(\bR)$:
\begin{equation}
  \label{eqn:annmult}
  A w(t)= a(t)w(t), \, t\in \bR.
\end{equation}
Explicit choices
for $a$ are discuss below.

With these choices, define
\begin{eqnarray}
  \label{eqn:edef}
  e[m,w;d] & = & \frac{1}{2}\|F[m]w-d\|^2/\|d\|^2;\\
  \label{eqn:gdef}
  g[w;d] & = & \frac{1}{2}\|Aw\|^2/\|d\|^2;\\
  \label{eqn:jdef}
  \Ja[m,w;d] & = & e[m,w;d] + \alpha^2g[w;d].
\end{eqnarray}

The main theoretical device used in the proofs of our main results on 
extended inversion is Variable Projection reduction of the penalty objective $\Ja$
(equation \ref{eqn:jdef}) to a function $\tJa$ of 
$m$ alone, by minimization over $w$:
\begin{equation}
  \label{eqn:redexp}
  \tJa[m;d] = \inf_w \Ja[m,w;d]
\end{equation}
A minimizer $w$ on the right-hand side of definition
\ref{eqn:redexp} must solve the {\em normal equation}
\begin{equation}
  \label{eqn:norm}
  (F[m]^TF[m]+\alpha^2A^TA)w= F[m]^Td, 
\end{equation}
(here and in the following, the superscript $T$ signifies adjoint in
the sense of domain and range Hilbert structures).

With $A$ of the form \ref{eqn:annmult}, $\tilde{J}_{\alpha}$ is explicitly
computable. First observe that apart from amplitude, $F[m]$ is
unitary: for $g \in D$,
\begin{equation}
\label{eqn:tran}
F[m]^T g (t) =
\left\{
  \begin{array}{c}
    \frac{1}{4\pi r}g\left(t+mr\right), \, t \in [t_{\rm min}-mr,
    t_{\rm max}-mr],\\
    0, \mbox{ else.}
  \end{array}
\right.
\end{equation}
so
\begin{equation}
  \label{eqn:unit}
  F[m]^TF[m] = \frac{1}{(4\pi r)^2}{\bf 1}_{[t_{\rm min}-mr,  
    t_{\rm max}-mr]}
\end{equation}
in which ${\bf 1}_{S}$ denotes
multiplication by the characteristic function of a measurable 
$S \subset \bR$.

Therefore the normal equation for the minimizer on the RHS of equation \ref{eqn:redexp} is
\begin{equation}
  \label{eqn:norm1}
  \left(\frac{1}{(4\pi r)^2} {\bf 1}_{[t_{\rm min}-mr,  
      t_{\rm max}-mr]} + \alpha^2 A^TA\right)w= F[m]^Td.
\end{equation}

With these choices, the normal equation \ref{eqn:norm1} becomes
\begin{equation}
\label{eqn:norm2}
\left(\frac{1}{(4\pi r)^2}  {\bf 1}_{[t_{\rm min}-mr,  
      t_{\rm max}-mr]} + \alpha^2a^2\right)w= F[m]^Td.
\end{equation}
Evidently, $F^TF$ (the first summand in equation \ref{eqn:norm2}) is
not coercive. If $A$ is to yield small output for input localized near
$t=0$, as it must to play its desired role in replacing the hard
support constraint, then $A^TA$ must not be coercive either, that is,
its spectrum must have zero as its infimum. Therefore an additional
condition on $A$ is required to enable application of the Lax-Milgram
constrruction to the normal
equation \ref{eqn:norm}. This condition is an hypothesis of 

\begin{proposition}
  \label{thm:norminvexp}
  Assume the conditions \ref{eqn:mod}, \ref{eqn:fullrec},
  \ref{eqn:annmult}. Also assume that $\lambda \in (0,\lambda_{\rm
    max}], \alpha > 0,$ and that  $C>0$ exists so that $a \in L^{\infty}(\bR)$
  mentioned in equation \ref{eqn:annmult} satisfies the condition
  \begin{equation}
    \label{eqn:abnd} 
    a \ge 0; \, a(t) \ge C\mbox{ for }|t| \ge \lambda_{\rm max}.
  \end{equation}
  Then
  \begin{itemize}
  \item[1. ]the normal operator $F[m]^TF[m] + \alpha^2A^TA$ is
    invertible for any $m \in M$, $\alpha > 0$;
  \item[2. ]the solution $\aw[m;d]\in W$ of the normal equation
    \ref{eqn:norm} is given by
    \begin{equation}
      \label{eqn:normsol}
      \aw[m;d](t) = \left\{
        \begin{array}{c}
          \left(\frac{1}{(4\pi r)^2} + \alpha^2
          a^2(t)\right)^{-1}\frac{1}{4 \pi r}d(t+mr), t \in [t_{\rm
          min}-mr, t_{\rm max}-mr];\\
          0, \mbox{ else;}
        \end{array}
      \right.
    \end{equation}
  \item[3. ]if in addition $d=F[m_*]w_*, w_* \in \lW$ is noise-free, as in equation
    \ref{eqn:defdatanonoise},
    \begin{equation}
      \label{eqn:solnonoise}
      \aw[m,d](t)= \left(1+ (4\pi r)^2\alpha^2 a(t)^2\right)^{-1}w_*\left(t+(m-m_*)r\right).
    \end{equation}
  \end{itemize}
\end{proposition}

\begin{proof}
  The idea here is that $F$ fails to be coercive just where $A$ is
  coercive, and vis-versa.
  \begin{itemize}
  \item[1. ]Note that thanks to \ref{eqn:zeroinc}, if $|t|\le
    \lambda \le \lambda_{\rm max}$, then ${\bf 1}_{[t_{\rm min}-mr,  
      t_{\rm max}-mr]}(t) = 1$, whereas if $|t|>\lambda$,
    then $a(t) \ge C$, whence
    \[
      \frac{1}{(4\pi r)^2}  {\bf 1}_{[t_{\rm min}-mr,  
        t_{\rm max}-mr]} + \alpha^2a^2  \ge \min\{(4\pi r)^2,
      \alpha^2\min\{1/(4\pi r)^2,C^2\} \}> 0.
    \]
    Therefore the normal operator is invertible under the stated
    conditions.

  \item[2. ]From the identity \ref{eqn:tran}.
    \[
      \mbox{supp }F[m]^Td \subset [t_{\rm min}-mr,t_{\rm max}-mr].
    \]
    Define $w_{\rm tmp}$ to be the right-hand side of equation \ref{eqn:normsol}. Then
    from the previous observation and identity \ref{eqn:tran},
    \[
      \mbox{supp }w_{\rm tmp} \subset [t_{\rm min}-mr,t_{\rm max}-mr].
    \]
    From the identity \ref{eqn:unit}, for any $w \in W$,
    \[
      t \in [t_{\rm min}-mr,t_{\rm max}-mr] \Rightarrow F[m]^TF[m]w(t)
      = \frac{1}{(4 \pi r)^2}w(t).
    \]
    It follows from this and the previous two observations that
    $w_{\rm tmp}$ solves the normal equation \ref{eqn:norm}, and
    therefore that $\aw[m;d]=w_{\rm tmp}$.

  \item[3. ]Follows by inserting the definition
    \ref{eqn:defdatanonoise} of $d$ in \ref{eqn:normsol} and
    rearranging.
  \end{itemize}
\end{proof}

\begin{theorem}
  \label{thm:norminv}
  Assume the condition \ref{eqn:fullrec}, $C>0$, and suppose that $A$ is
  given by equation \ref{eqn:annmult} for $a \in L^{\infty}(\bR)$
  satisfying condition \ref{eqn:abnd}.  Then
  \begin{itemize}
  \item[1. ]the reduced objective $\tJa$ is given by
    \begin{equation}
      \label{eqn:redexp1}
      \tJa[m;d] = \Ja[m,\aw[m;d];d],
    \end{equation}
    in which $\aw[m;d] \in W$ is the unique solution of the normal
    equation \ref{eqn:norm}.
  \item[2. ]The following are equivalent:
    \begin{itemize}
    \item[i. ]$(m,w) \in M \times W$ is a local minimizer of
      $\Ja[\cdot,\cdot;d]$, and
    \item[ii. ]$m$ is a local minimizer of $\tJa[\cdot;d]$ and
      $w=\aw[m;d]$.
    \end{itemize}
  \end{itemize}
\end{theorem}

\begin{proof} These conclusions follow immediately from Proposition
  \ref{thm:norminvexp}.
\end{proof}

If $\Ja[\cdot,\cdot;d]$ and $\tJa[\cdot;d]$
were differentiable, then ``local minimizer'' in the conclusion of the
preceding theorem could be replaced by ``stationary point''. However,
for the problem addressed in this paper, $\Ja[\cdot,\cdot;d]$ {\em is
  not} differentiable without added smoothness constraints on $w$,
whereas $\tJa[\cdot;d]$ {\em is} differentiable for proper choice of
penalty operator $A$.
This conclusion follows from properties of the modeling operator $F$
shared with many other inverse problems in wave propagation, as
explained in the Appendix. Here, I derive differentiability from explicit
expressions for $\tJa$ and its components.

\begin{proposition}
  \label{thm:epjgen}
  Assume the hypotheses of Proposition \ref{thm:norminvexp}. Then
  \begin{equation}
  \label{eqn:residnormgen}
  e[m,\aw[m,d];d] = \frac{1}{2\|d\|^2}\int_{t_{\rm min}}^{t_{\rm max}} \,dt\,(4\pi r \alpha a(t-mr))^4(1 +
  (4\pi r \alpha a(t-mr))^2)^{-2}d(t)^2
\end{equation}
\begin{equation}
  \label{eqn:anninormgen}
  p[m,\aw[m,d];d] = \frac{1}{2\|d\|^2}\int_{t_{\rm min}}^{t_{\rm max}} \,dt\,(4\pi r a(t-mr))^2(1 +
  (4\pi r \alpha a(t-mr))^2)^{-2}d(t)^2
\end{equation}

\begin{equation}
  \label{eqn:expjgen}
\tJa[m;d] = \frac{1}{2\|d\|^2}\int_{t_{\rm min}}^{t_{\rm max}}\,dt\,(4\pi r \alpha a(t-mr))^2(1+(4\pi r \alpha 
a(t-mr))^2)^{-1}d(t)^2. 
\end{equation}
\end{proposition}

\begin{proof}
  From equation \ref{eqn:normsol},
  \[
    F[m]\aw[m;d](t) = 
    \frac{1}{4 \pi r}\left(\frac{1}{(4\pi r)^2} + \alpha^2
      a^2(t-mr)\right)^{-1}\frac{1}{4 \pi r}d(t),
  \]
  so
  \[
    (F[m]\aw[m;d]-d)(t) = (1 + (4 \pi r\alpha
    a(t-mr))^2)^{-1}-1)d(t)
  \]
  \[
    = -(4 \pi r\alpha a(t-mr))^2(1 + (4 \pi r\alpha
    a(t-mr))^2)^{-1}d(t).
  \]
  Half the integral of the square of this data residual is
  $e[m,\aw[m;d],d]$, which proves identity \ref{eqn:residnormgen}.

  To compute $p[m,\aw[m;d],d]$, note that
  \[
    A\aw[m;d](t)=a(t) \left(\frac{1}{(4\pi r)^2} + \alpha^2
      a^2(t)\right)^{-1}\frac{1}{4 \pi r}d(t+mr)
  \]
  \[
    = 4\pi r a(t) (1 + (4\pi r \alpha a(t))^2)^{-1}d(t+mr)
  \]
  for $ t \in [t_{\rm min}-mr, t_{\rm max}-mr]$, so squaring,
  integrating, and changing integration variables $t \mapsto t-mr$
  gives the result \ref{eqn:anninormgen}

  That the VPM objective $\tJa$ is given by \ref{eqn:expjgen} follows from equations \ref{eqn:redexp1}, \ref{eqn:residnormgen}, and
  \ref{eqn:anninormgen}.
\end{proof}

\begin{theorem}
  \label{thm:diffobj}
  Suppose that in addition to the hypotheses of Theorem
  \ref{thm:norminv}, $a \in W^{1,\infty}_{\rm loc}(\bR)$, then $\tJa[\cdot;d]
  \in C^1(M)$.
\end{theorem}

\begin{proof}
Suppose first that $a \in C^1(\bR)$. Differentiation under the integral sign  
  yields the expression for its derivative:
\begin{equation}
  \label{eqn:dexpjgen}
  \frac{d}{dm}\tJa[m;d] = -\frac{(4 \pi r \alpha)^2}{\|d\|^2} \int_{t_{\rm min}}^{t_{\rm max}} \,dt \, 
  \left(a\frac{da}{dt}\right)(t-mr)(1+(4\pi r \alpha 
  a(t-mr))^2)^{-2}d(t)^2. 
\end{equation}
For $a \in W^{1,\infty}_{\rm loc}(\bR)$ a limiting argument shows that the
same expression gives the derivative of $\tJa$.
\end{proof}

It will be useful to record expressions for the various componenets of
$\tJa$ when the data is noise-free, that is, the context of
Proposition \ref{thm:norminvexp}, item 3:
\begin{corollary}
  \label{thm:epjnonoise}
  Assume the hypotheses of Proposition \ref{thm:norminvexp}, item
  3. Then noting that $\|d\| = \|w_*\|/(4 \pi r)$
\begin{equation}
  \label{eqn:residnorm}
  e[m,\aw[m,d];d] 
= \frac{\alpha^4}{2\|w_*\|^2}\int\,dt\,a(t-(m-m_*)r)^4(1+(4\pi r)^2 \alpha^2 
    a(t-(m-m_*)r)^2)^{-2}w_*(t)^2.
\end{equation}
\begin{equation}
  \label{eqn:anninorm}
  p[m,\aw[m,d];d] = \frac{(4\pi r)^2}{2\|w_*\|^2}\int \,dt\,  
  \frac{a(t-(m-m_*)r)^2}{(1+ (4\pi r)^2\alpha^2
    a(t-(m-m_*)r)^2)^{2}}w_*(t)^2.
\end{equation}
so 
\begin{equation}
\label{eqn:expjnonoise}
\tJa[m;d] = \frac{(4\pi r \alpha)^2}{2\|w_*\|^2}\int\,dt\,a(t-(m-m_*)r)^2(1+(4\pi r)^2 \alpha^2 
  a(t-(m-m_*)r)^2)^{-1}w_*(t)^2. 
\end{equation}
Finally, if $a \in W^{1,\infty}(\bR)$, then $\tJa[\cdot;d]$ is differentiable, and 
\begin{equation}
  \label{eqn:dexpjnonoise}
  \frac{d}{dm}\tJa[m;d] = -\frac{r (4\pi r \alpha)^2}{\|w_*\|^2} \int \,dt \, 
  \frac{\left(a\frac{da}{dt}\right)(t-(m-m_*)r)}{(1+(4\pi r)^2 \alpha^2 
  a(t-(m-m_*)r)^2)^{2}}w_*(t)^2. 
\end{equation}
\end{corollary}

Recall that the purpose of the penalty operator is to penalize energy away from 
$t=0$. A simple multiplier of class $W^{1,\infty}$ that accomplishes
this goal (and produces a differentiable reduced objective $\tJa$,
thanks to Theorem \ref{thm:diffobj}) is
\begin{equation}
  \label{eqn:ann}
  a(t) = \min(|t|, \tau). 
\end{equation}
for suitable $\tau>0$. In fact, 
the cutoff $\tau$ will be chosen large enough to be effectively inactive: 
hindsight suggests 
\begin{equation}
  \label{eqn:taudef}
  \tau = \max\{|t_{\rm min}-m_{\rm min}r|,|t_{\rm min}-m_{\rm max}r|, |t_{\rm max}-m_{\rm min}r|, |t_{\rm max}-m_{\rm max}r|\}. 
\end{equation}
Essentially this 
particular annihilator has been employed in earlier work on extended 
source inversion 
\cite[]{Plessix:00a,LuoSava:11,Warner:14,HuangSymes:SEG15a,Warner:16,HuangSymes:Geo17}.

\section{Stationary Points}

\begin{proposition}
  \label{thm:rampgood}
  Suppose that
  \begin{enumerate}
  \item $m_* \in M$;
  \item $0 < \mu \le \lambda$, and $w_* \in W_{\mu}$;
  \item $d_* = F[m_*]w_*$;
  \item $a(t)=\min\{|t|,\tau\}$ in the definition \ref{eqn:annmult},
    with $\tau$ given by equation \ref{eqn:taudef}; and
  \item $\alpha > 0$.
  \end{enumerate}
  Then for any $m \in M$, 
  \begin{equation}
    | (m - m_*)r| > \lambda  \Rightarrow  \left|\frac{d}{dm}\tJa[m;d_*]\right| >  
    \frac{r(4 \pi r \alpha)^2(\lambda-\mu)}{(1+(4\pi r\alpha)^2 
      (\lambda+\mu)^2)^{2}}  
    \label{eqn:gradbndnonoise}
  \end{equation}
\end{proposition}
\begin{proof}
  As observed before, $\mbox{supp }\aw[m;d_*] \subset [t_{\rm
    min}-mr,t_{\rm max}-mr]\subset [-\tau,\tau]$, with $\tau$ defined
  in \ref{eqn:taudef}. Therefore, $a(t) = |t|$, $a a'(t) = t$ in the
  support of the integrand on the RHS of equation
  \ref{eqn:dexpjnonoise}, which therefore 
  becomes (after change of integration variable)
  \begin{equation}
    \label{eqn:gradfinal}
    \frac{d}{dm}\tJa[m;d_*] = -\frac{r (4\pi r\alpha)^2}{\|w_*\|^2} \int \,dt \, 
  t(1+(4\pi r)^2 \alpha^2 
  t^2)^{-2}w_*(t+(m-m_*))^2.
  \end{equation}
  Recall that $w_*(t+(m-m_*)r)$
  vanishes if $|t+(m-m_*)r| > \lambda$. Therefore the integral on the
  RHS of equation \ref{eqn:gradfinal} can be re-written
  \[
    = -\frac{r(4 \pi r \alpha)^2}{\|w_*\|^2}\int_{-(m-m_*)r-\lambda}^{-(m-m_*)r+\lambda}
    \,dt\, t(1+(4\pi r)^2\alpha^2 t^2)^{-2}w_*\left(t+(m-m_*)r\right)^2
  \]
  Suppose that $\mu \le \lambda$ and $w_* \in W_{\mu}$. 
  If $m > m_*+\lambda/r$, then $t+(m-m_*)r \in \mbox{supp }w_*$
  implies $-\mu - \lambda < t < \mu-\lambda<0$, so 
  \[
    t(1+(4\pi r)^2\alpha^2 t^2)^{-2} < (\mu-\lambda)(1+(4\pi r)^2\alpha^2 (\mu+\lambda)^2)^{-2}<0
  \]
  in the support of the integrand in equation
  \ref{eqn:gradfinal}. Arguing similarly for $m<m_*-\lambda/r$, obtain
  a similar inequality, implying the conclusion \ref{eqn:gradbndnonoise}.
\end{proof}

\begin{theorem}
  \label{thm:rampreallygood}
  Suppose that
  \begin{enumerate}
  \item $m_* \in M$;
  \item $0 <  \lambda$, and $w_* \in \lW$;
  \item $d_* = F[m_*]w_*$;
  \item $a(t)=\min\{|t|,\tau\}$ in the definition \ref{eqn:annmult},
    with $\tau$ given by equation \ref{eqn:taudef}; 
  \item $\alpha > 0$; and
  \item$m \in M$ is a stationary point of $\tJa[\cdot;d_*]$.
  \end{enumerate}
  Then $|m-m_*| < \lambda /r$.
\end{theorem}

\begin{proof} Follows directly from Proposition \ref{thm:rampgood} by
  taking $\mu=\lambda$.
\end{proof}

The preceding theorem established that a proper choice of annihilator
leads to a reduced penalty objective all of whose stationary points
are within $O(\lambda)$ of the target slowness $m_*$, provided that
the data are noise-free in the sense of equation
\ref{eqn:defdatanonoise}. This result leaves open two questions:
\begin{itemize}
\item how does one use this reduced penalty minimization to produce
  a solution of the inverse problem as in problem statement
  \ref{eqn:probstat0}? 
\item how does one answer the same question for noisy data?
\end{itemize}

The next result answers the first question, in the case of noise-free data:
\begin{proposition}
  \label{thm:ipnonoisesuf}
  Suppose that $a$ is given by definition \ref{eqn:ann}, $\alpha$,
  $\mu \in (0,\lambda_{\rm max}]$,
  $d$ is given by
  \ref{eqn:defdatanonoise} with $w_* \in W_{\mu}$, and  $m$ is a stationary
  point of $\tJa[\cdot;d]$. Then $(m,\aw[m;d])$ is a
  solution of the inverse problem \ref{eqn:probstat0} for any $\lambda
  \ge 2\mu$ and
  \begin{equation}
    \label{eqn:estresidnorm}
    \epsilon \ge \frac{(8\pi r \mu \alpha)^2}{1 + (8\pi r \mu\alpha)^2}.
  \end{equation}
\end{proposition}

\begin{proof}
  From the assumption $w_* \in W_{\mu}$ and Theorem
  \ref{thm:rampreallygood}, $|(m-m_*)r|\le \mu$. From the
  identity \ref{eqn:solnonoise},
  $\mbox{supp }\aw[m;d] \subset
  [(m-m_*)r-\mu,(m-m_*)r+\mu] \subset
  [-2\mu,2\mu]$. Because of the support limitation, $a(t)=|t|$ in the
  interval of integration appearing in \ref{eqn:residnorm}, so
\[
  e[m,\aw[m,d];d] 
= 8 \pi^2 r^2 \alpha^4\int^{\mu}_{-\mu}\,dt\,\frac{|t-(m-m_*)r|^4}{(1+(4\pi r)^2 \alpha^2 
|t-(m-m_*)r|^2)^{2}}w_*(t)^2
\]
\[
  \frac{1}{2} (4\pi r \alpha)^4\int^{\mu}_{-\mu}\,dt\,\frac{|t-(m-m_*)r|^4}{(1+(4\pi r)^2 \alpha^2 
|t-(m-m_*)r|^2)^{2}}d(t+m_*r)^2
\]
\[
  \le \frac{1}{2} \|d\|^2  \left(\frac{(8\pi r \alpha \mu)^4}{(1+(8\pi
      r \alpha \mu)^2)^2}\right).
  \]
\end{proof}

The inequality \ref{eqn:estresidnorm} can be interpreted as a bound 
on $\alpha$, given $\epsilon$ and $\lambda$, for a
stationary point of $\tJa$ to yield a solution of the inverse
problem: one obtains a solution, provided that $\alpha$ is
sufficiently small. On the other hand, it is clear that $\alpha$
cannot be too large if stationary points of $\tJa$ are to yield
solutions: the integrand in \ref{eqn:residnorm} is increasing in
$\alpha$ for every $t$ and $m$, and the multiplier
\[
t \mapsto (4\pi r \alpha(t-(m-m_*)r))^4(1+(4\pi r)^2 \alpha^2 
|t-(m-m_*)r|^2)^{-2}
\]
tends monotonically to $1$ as $\alpha \rightarrow \infty$, uniformly
on the complement of any open interval containing
$t=(m-m_*)r$. Therefore
\begin{equation}
  \label{eqn:elimit}
  \lim_{\alpha \rightarrow \infty}e[m,\aw[m;d];d] =
  \frac{1}{2}\frac{1}{(4 \pi r)^2}\|w_*\|^2 = \frac{1}{2}\|d\|^2.
\end{equation}
Consequently, there exists $\alpha_{\rm max}(\epsilon,\lambda,d)$ so
that
\[
  e[m,\aw[m;d];d]  \le \frac{1}{2}\epsilon^2\|d\|^2
  \Rightarrow \alpha \le \alpha_{\rm max}(\epsilon,\lambda,d).
\]
The existence of this limiting penalty weight has been inferred
indirectly. \cite{FuSymes2017discrepancy} describe a constructive
algorithm for its approximation, which is used in the companion paper
\cite{SymesChenMinkoff:21} to dynamically adjust $\alpha$ in the
course of optimization of $\tJa$.

I turn now to the second issue identified above, the effect of
noise. Suppose that the data trace $d$ takes the form
\begin{equation}
  \label{eqn:defdatanoisy}
  d = F[m_*]w_* + n = d_*+n,
\end{equation}
with $m_* \in M, w_* \in W_{\mu}$, $0<\mu<\lambda$, and noise trace $n \in
D$. Since no support assumptions can be made about $n$, equation
\ref{eqn:normsol} implies that $\aw[m;d] \notin \lW$ for any values of
$\alpha$ and $\lambda$.  Therefore minimization of $\tJa$ cannot by itself yield a
solution of the inverse problem as defined in the problem statement
\ref{eqn:probstat0}. In this section, I explain how a solution may
nonetheless be constructed from a stationary point of $\tJa$.

First, examine the effect of additive noise on the estimation of the
slowness $m$. In expressing the result, use the dimensionless
relative data error
\begin{equation}
  \label{eqn:defeta}
  \eta = \frac{\|n\|}{\|d_*\|}. 
\end{equation}

\begin{proposition}
  \label{thm:mnoise}
  Assume the hypotheses of Proposition \ref{thm:rampgood}, and that $d$ is
  given by definition \ref{eqn:defdatanoisy}. Suppose that $m \in M$
  is a stationary point of $\tJa[\cdot;d]$, and that
  \begin{equation}
    \label{eqn:mnoisebnd}
    \eta(1+\eta) \le \frac{16}{3\sqrt{3}}\frac{4\pi r \alpha
      (\lambda-\mu)}{(1+(4\pi r\alpha(\lambda+\mu))^2)^2}
  \end{equation}
  Then
  \begin{equation}
    \label{eqn:mnoisebndfin}
    |m-m_*| \le \frac{\lambda}{r}
  \end{equation}.
\end{proposition}

\begin{proof}
  From equation \ref{eqn:dexpjgen}, $d \tJa / dm$ is the
value of a quadratic form in $d$ with (indefinite) symmetric operator
$B = $ multiplication by
\[
  b(t;m,\alpha)  = -\frac{(4 \pi r \alpha)^2 (t-mr)}{(1+(4\pi r \alpha (t-mr))^2)^{2}}
\]
Therefore
\begin{equation}
  \label{eqn:gradlip}
  \left|\frac{d}{dm}\tJa[m;d]-\frac{d}{dm}\tJa[m;d_*]\right| =
  |\langle (d+d_*),B(d-d_*)\rangle| \le \max_{t \in
    \bR}|b(t;m,\alpha)|\eta(1+\eta)\|d_*\|^2
\end{equation}
A straightforward calculation shows that
\[
  \max_{t \in \bR} b(t;m,\alpha) = \frac{3\sqrt{3}}{16} 4\pi r\alpha.
\]
For a stationary point $m$ of
$\tJa[\cdot;d]$, the inequality \ref{eqn:gradlip} implies
\[
  \left|\frac{d}{dm}\tJa[m;d_*]\right| \le \frac{3\sqrt{3}}{16} 4\pi
  r\alpha \eta(1+\eta)\|d_*\|^2
\]
On the other hand, the conclusion \ref{eqn:gradbndnonoise} of Proposition
\ref{thm:rampgood} implies that if also
\[
  \frac{3\sqrt{3}}{16} 4\pi r\alpha \eta(1+\eta)\|d_*\|^2 \le (4 \pi r
  \alpha)^2 \frac{\lambda-\mu}{(1+(4\pi r)^2\alpha^2
    (\lambda+\mu)^2)^{2}} \|d_*\|^2
\]
then $|m-m_*|\le \lambda/r$. Rearranging, obtain the conclusion.
\end{proof}

\begin{corollary}
  \label{thm:mnoisecor}
  Asumme the hypotheses of Proposition \ref{thm:mnoise}, in particular
  that $m$ is a stationary point of $\tJa[\cdot;d]$, $d=d_*+n$. Then
  \begin{equation}
    \label{eqn:mnoisecor}
    |m-m_*| \le \frac{\mu}{r} + \frac{\eta}{\alpha} \left(\frac{3\sqrt{3}(1+\eta)}{64\pi r^2}(1+(8\pi r \alpha
      \lambda_{\rm max})^2)^2\right)
  \end{equation}
\end{corollary}

\begin{proof} Assume that $\lambda$ is chosen to obtain equality in
  the condition \ref{eqn:mnoisebnd}, substitute the bound $2
  \lambda_{\rm max}$ for $\lambda + \mu$ in the denominator, solve for
  $\lambda$ and substitute in inequality \ref{eqn:mnoisebndfin}.
\end{proof}

\noindent {\bf Remark:} This bound suggests that error in the estimate
of $m$ due to data noise is a decreasing function of $\alpha$, at
least for small $\alpha$. This result is intuitively appealing, and is
supported by numerical evidence. It may be taken as a justification
for the discrepancy-based algorithm for adjusting $\alpha$ during
inversion \cite[]{FuSymes2017discrepancy}, used in the companion paper
\cite[]{SymesChenMinkoff:21}. 

\begin{theorem}
  \label{thm:mnoiseres}
  Asumme that
  \begin{itemize}
  \item[1. ] $\alpha, \mu> 0$,
  \item[2. ] $m_* \in M, w_* \in W_{\mu}$,
  \item[3. ] $d_* = F[m_*]w_*$,
  \item[4. ] $n \in D$ and $d = d_* + n$.
  \end{itemize}
  Set $\eta = \|n\|/\|d_*\|$. Assume that $\eta$ satisfies inequality \begin{equation}
  \label{eqn:mnoisecond}    
  \eta < \frac{\sqrt{5}-1}{2},
  \end{equation}
  and that $m$ is a stationary point of $\tJa[\cdot;d]$.
  Then
  \begin{equation}
    \label{eqn:mnoisesuff}
    |m-m_*| \le \left(1+\frac{2\eta(1+\eta)}{1-\eta(1+\eta)}\right)\frac{\mu}{r}.
  \end{equation}
\end{theorem}

\begin{proof} of Theorem \ref{thm:mnoiseres}:
  Write $\lambda = (1+\delta)\mu$, and $x=4 \pi r \alpha \mu$. Then
  the right-hand side of equation \ref{eqn:mnoisebnd} may be written as
  \begin{equation}
    \label{eqn:mnoisebndrev}
    \frac{16}{3\sqrt{3}}\frac{4\pi r \alpha
      (\lambda-\mu)}{(1+(4\pi r\alpha(\lambda+\mu))^2)^2} = D
    \frac{x}{(1+C^2 x^2)^2},
  \end{equation}
  where
  \[
    D=\frac{16}{3\sqrt{3}}\delta,\,C=2+\delta.
  \]
  The positive stationary point of the quantity on the right-hand side
  of \ref{eqn:mnoisebndrev} is a maximum, and occurs at
  $x=1/(\sqrt{3}C)$, that is
  \[
    4 \pi r \alpha \mu = \frac{1}{\sqrt{3}(2+\delta)}.
  \]
  Thus
  \[
    1+C^2x^2 = \frac{4}{3}
  \]
  hence the maximum value is
  \[
    \frac{D3\sqrt{3}}{16C} = \frac{\delta}{2+\delta}.
  \]
  This maximum value must be larger than the left hand side of inequality
  \ref{eqn:mnoisebnd}, that is,
  \[
    \eta(1+\eta) \le \frac{\delta}{2+\delta},
  \]
  in order that there be any solutions at all, but the right hand side
  is less than $1$. This observation establishes the necessity of
  hypothesis \ref{eqn:mnoisecond} of the
  theorem. Solving this above inequality for $\delta$ and unwinding
  the definitions, one finds that the right-hand side of inequality
  \ref{eqn:mnoisebnd} is bounded by \ref{eqn:mnoisesuff}, so appeal to
  Proposition \ref{thm:mnoise} finishes the proof.
\end{proof}

\noindent {\bf Remark.} That is, {\em with the choice of penalty
  multiplier $a$ given in equation \ref{eqn:ann}, and support radius
  $\mu$ of the ``noise-free'' wavelet $w_*$, $\Ja$ has no local
  minima with slownesses further than $(1+O(\eta))\mu/r)$ from the slowness used to
  generate the data}.

\noindent {\bf Remark.} The estimate $|m-m_*|r<\mu(1+O(\eta))$ for
local minima of $\Ja$ is sharp: it is possible to choose $w_* \in W_{\mu}$
so that $\mu - |m-m_*|r$ is as small as you like. In particular,
the ``exact'' or ``true'' slowness $m_*$ is not necessarily the only
slowness component of a local minimizer, or even the slowness
component of any local minimizer, and in particular
is not (necessarily) the slowness component of a global minimizer of $\Ja$.

\noindent {\bf Remark.} No similar bound could hold for much larger
noise levels than specified in condition \ref{eqn:mnoisecond}, the
right-hand side of which is a bit larger than 0.6. For example, if the noise is
the predicted data for the same wavelet $w_*$ with a substantially different
slowness $m_{\flat}$, that is, $n=F[m_{\flat}]w_*$, then a simple
symmetry argument shows that if there is a local minimizer of $\Ja[\cdot,\cdot;d_*+n]$. with
slowness near $m_*$,
there must also be a minimizer with slowness near $m_{\flat}$.
so that the difference with $m_*$ is not constrained at
all by the assumed support radius of $w_*$. So for this example with 100\% noise, no
bound of the type given by conclusion 2 could possibly hold. The
companion paper \cite[]{SymesChenMinkoff:21}
illustrates this phenomenon numerically.

\noindent{\bf Remark.} Note that $\alpha$ plays no role in the
conclusions of this theorem. It is only required that $\alpha >0$.

\noindent {\bf Remark:} I emphasize that Theorem \ref{thm:mnoiseres} states {\em sufficient} conditions for a bound on
the slowness error $|m-m_*|$ in terms of the relative data noise level $\eta$,
giving an additional ``fudge factor'' beyond the support size $\mu$
of the noise-free wavelet $w_*$ for an interval within which the slowness error is
guaranteed to lie.

Conclusion 1 in Theorem \ref{thm:mnoiseres} constrains the range of
noise level to which these results apply to a bit more than 60\%. That
is, the bound given by conclusion 2 is useful only for small noise. In
the limit as $\eta \rightarrow 0$, conclusion 2 becomes
$\lambda/\mu \gtrsim 1 + 2\eta$, that is, the ``fudge factor'' beyond
the noise-free bound is approximatly twice the noise level.

On the other hand, stronger bounds than given by Theorem
\ref{thm:mnoiseres} are possible, given additional constraints on the
noise $n$. A natural example is uniformly distributed random noise,
filtered to have the same spectrum as the source. The expression
\ref{eqn:dexpjgen} implies that the interaction of noise $n$ and
signal $d_*$ in the derivative of $\tJa$ is local, so that the
coefficient of $\eta$ on the left-hand side of inequality
\ref{eqn:mnoisebnd} is effectively much less that 1, resulting in a
larger range of allowable $\eta$. While I will not formulate such a
result, one of the numerical examples in the companion paper
\cite[]{SymesChenMinkoff:21} suggests its feasibility.

Unless the data is noise-free, there is no reason to suppose that the
estimated wavelet $\aw[m;d]$ (Theorem \ref{thm:norminv}) will lie in
$\lW$, unless the support of the noise $n$ is restricted. In order
to construct a solution of the inverse problem \ref{eqn:probstat0}, 
project $\aw[m;d]$ onto $\lW$. For sufficiently large $\lambda,
\epsilon$, the result is a solution of the inverse problem:

\begin{theorem}
  \label{thm:ipnoisesuf}
  Assume the hypotheses of Theorem \ref{thm:mnoiseres}, and that
  inequality \ref{eqn:mnoisecond} holds, and $\mu \in
  (0,\lambda_{\rm max}]$. Then the pair
  \[
    (m,{\bf 1}_{[-\lambda,\lambda]}\aw[m,d])
  \]
  solves the inverse problem as stated in \ref{eqn:probstat0} if
  \begin{equation}
    \label{eqn:ipnoiselam}
    \left( 2+\frac{2\eta(1+\eta)}{1-\eta(1+\eta)}\right)\mu \le \lambda
    \le \lambda_{\rm max}, 
  \end{equation}
  and
  \begin{equation}
    \label{eqn:ipnoiseeps}
    \epsilon \ge \frac{(8 \pi r \alpha\lambda)^2}{1 + (8 \pi r \alpha\lambda)^2}+\eta. 
  \end{equation}    
\end{theorem}

\begin{proof} of Theorem \ref{thm:ipnoisesuf}:
From Theorem \ref{thm:norminv}, 
\[
  \aw[m;d](t) = (1+ (4 \pi r \alpha t)^2)^{-1}(w_*(t+(m-m_*)) + 4\pi r 
  n(t+mr)) 
\]
\[
  = \aw[m;d_*] + (1+ (4 \pi r \alpha t)^2)^{-1}4\pi r n(t+mr)
\]

From Theorem \ref{thm:mnoiseres},
\[
  |m-m_*| \le \left(1+\frac{2\eta(1+\eta)}{1-\eta(1+\eta)}\right)\frac{\mu}{r} 
\]
which from assumption \ref{eqn:ipnoiselam} is
\[
  =\left(2+\frac{2\eta(1+\eta)}{1-\eta(1+\eta)}\right)\frac{\mu}{r}
  -\frac{\mu}{r} \le \left(\frac{\lambda}{\mu}-1\right)\frac{\mu}{r}
\]
That is,
\[
  |m-m_*|r \le \lambda-\mu.
\]
From Theorem \ref{thm:norminv},
\[
  \mbox{supp }\aw[m;d_*] \subset [-\mu-(m-m_*)r, \mu-(m-m_*)r]
  \subset [-\lambda,\lambda]
\]
so
\[
  E_{\lambda}\aw[m;d](t) = \aw[m;d_*](t) + E_{\lambda}(1+ (4 \pi r \alpha
  t)^2)^{-1}4\pi r n(t+mr)
\]
From the definition of $F[m]$, for any $t_1<t_2$, $w \in W$,  
\[
  F[m]{\bf 1}_{[t_1.t_2]}w = {\bf 1}_{[t_1+mr,t_2+mr]}F[m]w  
\]
Thus the data residual after projection is
\[
  F[m]E_{\lambda}\aw[m,d](t) -d(t) = F[m]\aw[m,d_*](t) -d_*(t)  
\]
\[
  + {\bf 1}_{[-\lambda+mr,\lambda+mr]}(4 \pi r \alpha (t-mr))^2 (1+ (4 \pi
  r \alpha (t-mr))^2)^{-1} n(t)
\]
\[
  -(1- {\bf 1}_{[-\lambda+mr,\lambda+mr]})n(t)
\]
From \ref{eqn:residnorm} and the bound on $m-m_*$,
\[
  \| F[m]E_{\lambda}\aw[m,d_*] -d_*\|^2 = (4 \pi r \alpha)^4
  \int_{-\lambda+(m-m_*)r}^{\lambda+(m-m_*)r}\,dt\, (t-(m-m_*)r)^4
\]
\[
  \times (1+(4\pi r \alpha)^2 
  (t-(m-m_*)r)^2)^{-2}d_*(t)^2
\]
\[
  \le (4 \pi r \alpha \lambda)^4\|d_*\|^2
\]
Similarly, the norm squared of the sum of the last two terms is 
\[
  \le (4 \pi r \alpha \lambda)^2 \|{\bf 1}_{[-\lambda+mr,\lambda+mr]}
  n\|^2 + \|(1 - {\bf 1}_{[-\lambda+mr,\lambda+mr]}n\|^2
\]
Without additional hypotheses to outlaw the accumulation of $n$ near
$t=mr$, all that can be said is that this is
\[
  \le \max \{(4 \pi r \alpha \lambda)^2, 1\} \|n\|^2
\]
Putting this all together,
\[
  \|F[m]E_{\lambda}\aw[m,d]-d\| \le (4 \pi r \alpha
  \lambda)^2\|d_*\| + \max \{4\pi r \alpha \lambda, 1\}\|n\|
\]
\[
  \le ((4 \pi r \alpha \lambda)^2 +\max \{4\pi r \alpha \lambda, 1\}
  \eta)\|d_*\|.
\]
If the right-hand side is to be less than $\|d_*\|$ as required by the
definition \ref{eqn:probstat0} of the inverse problem, then
necessarily $4\pi r\alpha \lambda < 1$, so the right hand side in the
preceding inequality is bounded by the right hand side of assumption
\ref{eqn:ipnoiseeps} of the theorem. Therefore this assumption implies
that the relative residual is $\le \epsilon$.
\end{proof}

\noindent {\bf Remark:} Note that the sufficient condition \ref{eqn:ipnoiselam} for
$\lambda$ is independent of $\alpha$. It follows that for any choice
of $\lambda$ consistent with this bound, $(m,{\bf
  1}_{[-\lambda,\lambda]}\aw[m,d])$ is a solution of the inverse
problem for any $\epsilon > 0$ provided that $\alpha$ is chosen
sufficiently small ($O(\sqrt{\epsilon})$).

\bibliographystyle{seg}
\bibliography{../../bib/masterref}

\section*{APPENDIX: ABSTRACT STRUCTURE OF THE GRADIENT}

The inverse problem studied here is far too simple to have any direct
use in applications. Its simplicity allows a rather complete account of
its properties, as the reader has seen in the preceding pages. However
this discussion has touched on features found in more complex and
prototypical problems:
\begin{itemize}
\item Nested optimization (variable projection method,
  \cite{GolubPereyra:03}) based on a model decomposition into inner
  and outer variables seems to be very important: in some cases, such
  as the simple problem presented here, the decomposition is obvious
  (inner variable = wavelet, outer variable = slowness), in others
  less so \cite[]{geoprosp:2008,Terentyev:thesis}.
\item The extended modeling operator must be surjective, or at least
  have a dense range, {\em for each value of the outer variable}, so
  that data may be fit well even for a poor initial guess of the outer
  variable. for the extended inversion approach to be successful.
\item Because of this data-fitting assumption, a straightforward
  algorithm for scaling the quadratic penalty, based on the
  Discrepancy Principle, is available - see
  \cite{FuSymes2017discrepancy,Symes:21a,SymesChenMinkoff:21}. This
  scaling varies dynamically during iterative optimization, in effect
  changing the objective function sporadically as the iteration
  converges.
\item The derivative of the extended modeling operator with respect to
  the outer variable is well-approximated by the composition of the
  extended modeling operator itself and a pseudodifferential operator
  of order 1 \cite[]{Symes:IPTA14,tenKroode:IPTA14}.
\item The choice of the penalty operator, controlling the extended
  degrees of freedom, is critical: in order to produce an objective
  immune from cycle-skipping, this penalty operator must be
  (pseudo-) differential of order zero \cite[]{StolkSymes:03}. A smooth
  multiplier $a$, as in equation \ref{eqn:annmult}, is a special case.
\end{itemize}

The expression \ref{eqn:dexpjgen} for the derivative of the reduced
objective $\tJa$ is the result of elementary manipulations, based on the
explicit expression \ref{eqn:mod} for the modeling operator $F$. In
this appendix I give an alternative derivation that generalizes to
extended inversion formulations for much less constrained physics. I
will point out the additional reasoning necessary to reach similar
conclusions in these more complex instances of extended inversion, as
presented for example in \cite{StolkSymes:03,StolkDeHoopSymes:09,Symes:IPTA14,tenKroode:IPTA14,HuangSymes:Geo17,HuangSymes:Geo18a,HuangSymes:Geo18b,HuangNammourSymesDollizal:SEG19}.

Recall that $\aw[m;d]$ is the solution of the normal equation
\ref{eqn:norm}. The reduced objective $\tJa$ is given by
\[
  \tJa[m;d] = \Ja[m,\aw[m;d];d]
\]
\[
  = \frac{1}{2}(\|F[m]\aw[m;d]-d\|^2 + \alpha^2\|Aw[m;d]\|^2)
\]
(equation \ref{eqn:redexp1})
\begin{equation}
  \label{eqn:redsimple}
  =\frac{1}{2}(\|d\|^2 - \langle d,F[m]\aw[m;d]\rangle),
\end{equation}
after a little algebra.

As mentioned above, $F:M \times W \rightarrow D$ is {\em not}
differentiable. Neither is $\aw: M \times D \rightarrow W$, as follows
immediately from the identity \ref{eqn:normsol}. However $F\aw: M
\times D \rightarrow D$ is differentiable, hence so is $\tJa: M \times
D \rightarrow \bR$ thanks
to the identity \ref{eqn:redsimple},
under the conditions on the multiplier $a$ identified in Theorem
\ref{thm:diffobj}.

For the model problem studied in the body of this paper, a simple
argument justifies this conclusion.  To sideline some technical
details, assume that
$a \in C^{\infty}(\bR)\cap L^{\infty}(\bR)$ and $a(t)=t$ for $|t|<\tau$ and $|a(t)| \ge
\tau$ for $|t| \ge \tau$, with $\tau$
defined in equation \ref{eqn:taudef}. Then
Proposition \ref{thm:norminvexp}, item [1], implies that
\begin{equation}
  \label{eqn:faw}
  F[m]\aw[m;d] = F[m](F[m]^TF[m] + \alpha^2 A^TA)^{-1}F[m]^T d.
\end{equation}
The normal operator and its inverse are multiplication
operators, whereas $F$ and $F^T$ are scaled shift operators, inverse
to each other except for scale. From this it is easy to see that the
operator on the RHS of equation \ref{eqn:faw} is multiplication by a
smooth function, with its arguments shifted by $mr$. Such an operator
is smooth in $m$, hence so is $F\aw$.

Since $\tJa$ and its gradient depend smoothly on their arguments, it
is possible to derive an alternate expression assuming that $d \in
H^1(\bR)$, then extend it by continuity to $d \in L^2(\bR)$. Note that
if $w \in H^1(\bR)$, then
  differentiable, and 
\begin{equation}
\label{eqn:deriv}
(D(F[m]w)\delta m)(t) = F[m](Q[m]\delta m)w (t), 
\end{equation}
where 
\begin{equation}
\label{eqn:defq}
(Q[m]\delta m)w = -r\delta m \frac{dw}{dt}. 
\end{equation}

That is, $DF[m]\delta m$ factors into $F[m]$ following $Q[m]\delta m$,
where the latter a skew-adjoint differential operator of order 1,
depending linearly on $\delta m$.

Assume that $d \in H^1(\bR)$. From equation \ref{eqn:normsol}, it
follows that $w[m;d] \in H^1(\bR)$ and moreover that $m \mapsto
\aw[m;d]$ is differentiable as a map from $\bR^+$ to
$L^2(\bR)$. From equation \ref{eqn:redexp1},
\[
  \tJa[m;d] = \frac{1}{2}(\|F[m]w[m;d]-d\|^2 + \alpha^2 \|Aw[m;d]\|),
\]
whence $m \mapsto \tJa[m;d]$ is differentiable. A standard calculation
invoking the normal equation \ref{eqn:norm} shows that
\[
  D\tJa[m;d]\delta m = \langle D(F[m]\aw[m;d])\delta m, F[m]\aw[m;d]-d\rangle
\]
\[
  = \langle F[m](Q[m]\delta m)\aw[m;d], F[m]\aw[m;d]-d\rangle
\]
\[
  = \langle (Q[m]\delta m)\aw[m;d], F[m]^T(F[m]\aw[m;d]-d)\rangle
\]
\[
  = -\alpha^2\langle (Q[m]\delta m)\aw[m;d],A^TA\aw[m;d] \rangle
\]
\[
  = \alpha^2 \langle \aw[m;d], (Q[m]\delta m)A^TA \aw[m;d]\rangle
\]
(using antisymmetry of $Q$)
\[
  = \alpha^2 \langle \aw[m;d],[ (Q[m]\delta m), A^TA]\aw[m;d] \rangle
  + \alpha^2 \langle (Q[m]\delta m)\aw[m;d],A^TA\aw[m;d] \rangle
\]
Rearranging,
\begin{equation}
  \label{eqn:djq}
  d\tJa[m;d]\delta m = \frac{1}{2}\alpha^2 \langle \aw[m;d],[
  (Q[m]\delta m), A^TA]\aw[m;d] \rangle.
\end{equation}
Since $Q$ is a differential operator of order 1, and $A^TA$ an
operator of order 0, the commutator has order 0. That is, the RHS of
equation \ref{eqn:djq} defines a continuous quadratic form in $d$ with
respect to the $L^2$ norm. As shown above, the same is true of the
left hand side. Therefore their extensions by continuity to $d \in
L^2(\bR)$ are the same, and the identity \ref{eqn:djq} holds for $d
\in L^2(\bR)$.

Making the identity \ref{eqn:djq} explicit by means of equation
\ref{eqn:defq} and $A^TA = t^2$, substituting the expression
\ref{eqn:normsol} for $\aw[m;d]$, and rearranging: one obtains
precisely the expression \ref{eqn:dexpjgen}.

Similar reasoning can be applied to some of the other extended
inversion settings mentioned at the beginning of this
appendix. Provided that the outer parameter $m$ consists of (or
parametrizes) smooth coefficients in the principal part of the
governing system of wave equations, $F$ is a microlocally elliptic
Fourier Integral Operator (FIO) \cite[]{Dui:95}. The canonical
relation of $F$ takes on the role of the shift operator
$t \rightarrow t-mr$, and has the properties necessary to conclude
that the composition of $F$ and $F^T$ in both orders is
pseudodifferential, at least in an open conic subset of the cotangent
bundle. If $A^TA$ is pseudodifferential, then so is the normal
operator, and it is microlocally elliptic, so positive-definite when
restricted to a suitable subspace of the inner parameter $w$, which
appears in the hyperbolic system as components of a right-hand side or
boundary condition. The inverse in the RHS of equation \ref{eqn:faw}
must be replaced by a microlocal parametrix. Thanks to Egorov's
Theorem \cite[]{Tay:81}, a special case of the rules for composing
FIOs, the operator on the RHS of equation \ref{eqn:faw} is a
pseudodifferential operator whose symbol is algebraic in geometric
optics quantities, hence depends smoothly on the coefficients in the
system of wave equations, and therefore on $m$, as does $F\aw$.

Microlocal ellipticity of $F$ implies an analogue of relation
\ref{eqn:deriv}, in which $Q$ is an essentially skew-symmetric
pseudodifferential operator of order 1. Since the derivation 
the gradient identity \ref{eqn:djq} is algebraic, it holds
approximately in these more general cases, with lower-order (smoother)
error terms. See 
\cite{Symes:Madrid,tenKroode:IPTA14,Symes:IPTA14} for details in
several specific cases, and for a similar calculation of the Hessian.
For these examples, the principal symbol of $Q$ can be computed
explicitly, as a function of the geodesic distance (traveltime)
function. Consequently, the analog of $\tJa$ is tangent to second
order to the Hessian of an objective formulated in terms of
traveltime. This relation explains the ability of the variant of
extended inversion studied in the cited references to recover
kinematically accurate velocity fields. Of course, that is what is
shown in detail in this paper for the very simple model problem
studied here.
\end{document}